\documentclass[journal]{IEEEtran}
\usepackage{tikz}
\usepackage[cmex10,intlimits]{amsmath}
\usepackage{amsfonts}
\usepackage{amssymb}
\usepackage{booktabs}

\graphicspath{figures/}

\usepackage[colorlinks]{hyperref}
\hypersetup{
	colorlinks=true, 
	linkcolor=blue!70!black, 
	citecolor=red!70!black, 
	filecolor=blue!70!black, 
	urlcolor=blue!70!black, 
}

\RequirePackage{xcolor}

\RequirePackage[cmex10,intlimits]{amsmath}
\RequirePackage{amsfonts}
\RequirePackage{amssymb}

\providecommand{\J}{\ensuremath{\mathrm{j}}}    

\providecommand{\Quot}[1]{``{#1}"}              
\providecommand{\V}[1]{\boldsymbol{#1}}         
\providecommand{\M}[1]{\mathbf{#1}}             
\providecommand{\T}[1]{\mathrm{#1}}             
\providecommand{\OP}[1]{{\mathcal{#1}}}         

\providecommand{\rv}{\ensuremath{\V{r}}}
\providecommand{\herm}{\mathrm{H}} 
\providecommand{\trans}{\mathrm{T}}

\providecommand{\srcRegion}{\ensuremath{\varOmega}}
\providecommand{\basisFcn}{\V{\psi}}

\providecommand{\Iv}{\ensuremath{\M{I}}}
\providecommand{\Vv}{\ensuremath{\M{V}}}

\providecommand{\Xm}{\ensuremath{\M{X}}}
\providecommand{\Ym}{\ensuremath{\M{Y}}}
\providecommand{\Zm}{\ensuremath{\M{Z}}}
\providecommand{\Zmvac}{\ensuremath{\M{Z}_0}}
\providecommand{\Rmvac}{\ensuremath{\M{R}_0}}
\providecommand{\Xmvac}{\ensuremath{\M{X}_0}}

\providecommand{\Wm}{\ensuremath{\M{W}}}
\providecommand{\Cm}{\ensuremath{\M{C}}}

\newcommand{\ie}{\textit{i}.\textit{e}.{}} 
\newcommand{\eg}{\textit{e}.\textit{g}.{}}
\newcommand{\cf}{\textit{cf}.{}}

\providecommand{\bbox}{\ensuremath{\varOmega_0}}
\providecommand{\gene}{\ensuremath{\M{g}}}

\providecommand{\setGi}{\ensuremath{\OP{E}}}
\providecommand{\setGo}{\ensuremath{\OP{B}}}
\providecommand{\setFixed}{\ensuremath{\OP{P}}}
\providecommand{\setAdd}{\ensuremath{\OP{A}}}
\providecommand{\setRem}{\ensuremath{\OP{R}}}

\providecommand{\Nopt}{\ensuremath{B}}
\providecommand{\Ndof}{\ensuremath{N}}
\providecommand{\Nags}{\ensuremath{N_{\mathrm{A}}}}

\RequirePackage[nolist]{acronym}
\newacro{MoM}[MoM]{method of moments}
\newacro{FEM}[FEM]{finite element method}
\newacro{FDTD}[FDTD]{finite-difference time-domain method}
\newacro{EFIE}[EFIE]{electric field integral equation}
\newacro{RWG}[RWG]{Rao-Wilton-Glisson}
\newacro{EM}[EM]{electromagnetic}
\newacro{DOF}[DOF]{\mbox{degrees-of-freedom}}

\begin{document}
\title{Optimal Inverse Design Based on Memetic Algorithms -- Part 1: Theory and Implementation}
\author{Miloslav~Capek, \IEEEmembership{Senior Member, IEEE}, Lukas~Jelinek, Petr~Kadlec, and Mats~Gustafsson, \IEEEmembership{Senior Member, IEEE}
\thanks{Manuscript received \today; revised \today. This work was supported by the Czech Science Foundation under project~\mbox{No.~21-19025M}.}
\thanks{M. Capek and L. Jelinek are with the Czech Technical University in Prague, Prague, Czech Republic (e-mails: \{miloslav.capek; lukas.jelinek\}@fel.cvut.cz).}
\thanks{M. Gustafsson is with Lund University, Lund, Sweden (e-mail: mats.gustafsson@eit.lth.se).}
\thanks{P. Kadlec is with the Brno University of Technology, Brno, Czech Republic (e-mail: kadlecp@vut.cz).}
\thanks{Color versions of one or more of the figures in this paper are
available online at http://ieeexplore.ieee.org.}
}

\maketitle

\begin{abstract}
A memetic framework for optimal inverse design is proposed by combining a local gradient-based procedure and a robust global scheme. The procedure is based on method-of-moments matrices and does not demand full inversion of a system matrix. Fundamental bounds are evaluated for all optimized metrics in the same manner, providing natural stopping criteria and quality measures for realized devices. Compared to density-based topology optimization, the proposed routine does not require filtering or thresholding. Compared to commonly used heuristics, the technique is significantly faster, still preserving a high level of versatility and robustness. This is a two-part paper in which the first part is devoted to the theoretical background and properties, and the second part applies the method to examples of varying complexity.
\end{abstract}


\section{Introduction}

\IEEEPARstart{I}{nverse} design is a long-lasting subject of study in electromagnetism. This is particularly true for antenna design~\cite{RahmatMichielssen_ElectromagneticOptimizationByGenetirAlgorithms, Haupt_Werner_GeneticAlgorithmsInEM}, the art of crafting and shaping the material and \ac{EM} sources to modify electric current paths so that they radiate effectively. Inverse design imposes two essential questions. The first is how good can, in principle, the design be. The second is how an optimal design would look like. These two questions are of different complexity but should be treated together.

The question regarding principal -- fundamental -- bounds was addressed by many researchers for various antenna metrics from bandwidth~\cite{Chu_PhysicalLimitationsOfOmniDirectAntennas} to antenna gain~\cite{UzsokySolymar_TheoryOfSuperDirectiveLinearArrays} to radiation efficiency~\cite{Harrington_EffectsOfAntennaSizeOnGainBWandEfficiency}. All works reflect the impact of electrically small design regions on all performance parameters~\cite{Fujimoto_Morishita_ModernSmallAntennas, VolakisChenFujimoto_SmallAntennas}. Yet, tightness of the bounds is known only in rare cases~\cite{Best2005b, Sievenpiper_ExpretimentalValidationOfPerformanceLimitsAndDesignGuidelines, Kildal+etal2017,GustafssonSohlKristensson_IllustrationsOfNewPhysicalBoundOnLinearlyPolAntennas}. This becomes problematic when computational inverse design techniques are used. When to stop them? Shall we rerun them? Clearly, the fundamental bounds and the way to approach them are closely interconnected problems.

The first step has been done in the field of fundamental bounds. Current-density-based bounds were proposed recently~\cite{Gustafsson_OptimalAntennaCurrentsForQsuperdirectivityAndRP} and formulated for many parameters~\cite{JelinekCapek_OptimalCurrentsOnArbitrarilyShapedSurfaces, GustafssonCapekSchab_TradeOffBetweenAntennaEfficiencyAndQfactor, GustafssonCapek_MaximumGainEffAreaAndDirectivity}. This class of bounds considers a specific design region and materials. The fundamental bound is associated with an optimal current density solution which is interpreted as being imposed on a hypothetical antenna structure lying in that region. Unfortunately, thanks to the relaxed constraints, the optimal current cannot be directly realized and, instead, an inverse design procedure has to be adopted to approach the performance predicted by the fundamental bound.

To address the second question, various procedures~\cite{HaslingerMakinen_IntroductionToShapeOptimization} have been proposed and met with mixed success, ranging from simple parametric sweeps to versatile heuristic algorithms~\cite{Simon_EvolutionaryOptimizationAlgorithms} to density-based topology optimization~\cite{BendsoeSigmund_TopologyOptimization}. All these methods are capable of delivering promising design candidates. Nevertheless, they embody serious practical weaknesses, mainly due to the NP-hard problems faced~\cite{Leeuwen_AlgorithmsAndComplexityA} stemming from the fact that the optimization is of combinatorial nature~\cite{Nemhauser_etal_IntegerAndCombinatorialOptimization}. Many improvements are known in the literature, \eg{}~\cite{Gregory_etal_FastOptimizationOfEMDesignProblemsUsingCovarianceMatrixAdaptionEvolutionaryStrategy, AldhafeeriRamatSamii_BrainStormOptimForEM, LiGuo_GrayWolfOptimForAntennaOptimDesign}, including memetic algorithms~\cite{Hart_etal_RecentAdvancesInMemeticAlgorithms}. However, they all are facing the no-free-lunch theorem~\cite{WolpertMacready_NoFreeLunchTheoremsForOptimization} stating that there is no unique technique suitable to all optimization problems.

A couple of assumptions are considered in this paper to propose an efficient optimization technique compatible with fundamental bounds. The underlying full-wave numerical method is the method of moments~\cite{Harrington_FieldComputationByMoM} with the possibility to use its matrices in the optimizer. Our attempt is to solve the original combinatorial problem~\cite{Ohsaki_OptimizationOfFiniteDimensionalStructures}, where an unknown shape is described with a characteristic function. The algorithm is a memetic combination of a local gradient method searching for local extrema and a global method maintaining diversity. Finally, the procedure should be scalable and fast. The last property is fulfilled by employing inversion-free evaluation based on block inversion of a matrix. 

Inversion-free evaluation was already proposed in~\cite{1989_Kastner_OnMatrixPartitioning} and later extended in~\cite{Capeketal_ShapeSynthesisBasedOnTopologySensitivity} proposing so-called topology sensitivity and further broadened in~\cite{Capeketal_InversionFreeEvaluationOfNearestNeighborsInMoM} by the possibility of structure growth. Several recent works have also been inspired by this approach~\cite{Jiang_etal_PixelAntenna_2022, WangHum_APS_BinaryTopologyModel}, sharing, however, the restrictions of the original approach, \ie{}, only removal is allowed, and only very small matrices can be optimized within single-purpose optimizer. A similar idea is used in~\cite{Budhu_Grbic_OptimizationGD_Woodbury} for the synthesis of metasurfaces or in~\cite{Zhang_etal_ABenchmarkTestSuiteForAntennaSparOpt} for benchmark testing suite for antenna S-parameter optimization.

To significantly boost the available number of unknowns necessary for effective antenna design, the techniques from~\cite{Capeketal_ShapeSynthesisBasedOnTopologySensitivity} and~\cite{Capeketal_InversionFreeEvaluationOfNearestNeighborsInMoM} are combined, and advanced version of pixeling procedure~\cite{RahmatSamii_Kovitz_Rajagopalan-NatureInspiredOptimizationTechniques, RahmatMichielssen_ElectromagneticOptimizationByGenetirAlgorithms, Haupt_Werner_GeneticAlgorithmsInEM} is added as an extra global step. Compared to the innovative local step, the global heuristic step is not the core of this technique and does not provide unique features. However, it provides additional diversity and movement within the solution space, strictly performed through the local minima. This is a unique feature of this work, not available elsewhere.

This paper is the first part of a two-part paper and is organized as follows. The optimization framework is introduced in Section~\ref{sec:ShapeOpt} consisting of fundamental bound evaluation and topology optimization. The local inversion-free update is introduced in Section~\ref{sec:Reanalysis}. The memetic combination of the local and global step (here implemented as a heuristic algorithm) is described in Section~\ref{sec:memetics}. The first part is concluded in Section~\ref{sec:concl}. The paper is accompanied by five appendices, showing mathematical details. In addition, a code implementing the local step is published as supplementary material~\cite{TSGA22_SuplMat}. In Part~2, the method's performance is demonstrated on several examples of varying complexity~\cite{2021_capeketal_TSGAmemetics_Part2}. All results are compared to the fundamental bounds based on convex optimization over current density.

\section{Optimal inverse design}
\label{sec:ShapeOpt}

\begin{figure}
\centering
\includegraphics[]{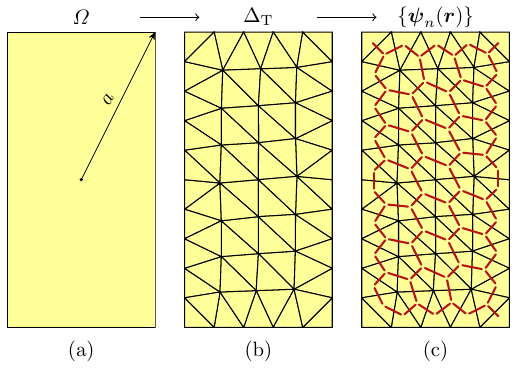}
\caption{(a) Region~$\varOmega$, here a rectangle with aspect ratio~$1:2$ circumscribed by a sphere of radius~$a$, (b) its discretization into a set of triangles~$\Delta_\T{T}$, and (c) edge basis function representation. The basis functions~$\left\{\basisFcn_n(\rv)\right\}$, $n \in \left\{1,\cdots,\Ndof\right\}$ are associated with the edges between two triangles as indicated by the red lines. Here, they are treated both as physical \acf{DOF} and optimization variables.}
\label{fig1}
\end{figure}

Without loss of generality, let us consider optimization of a surface structure. For this purpose, imagine region~$\varOmega$ available for the inverse design of a device, see Fig.~\ref{fig1}a. The numerical method used in this work is \ac{MoM} which requires to discretize the entire region~$\varOmega$, where the unknown polarization currents reside, see Fig.~\ref{fig1}b. In order to express the problem in a tractable algebraic way~\cite{Harrington_FieldComputationByMoM}, $N$ overlapping Rao-Wilton-Glisson basis functions~$\basisFcn_n(\rv)$ (defined over every doublet of adjacent triangles)~\cite{RaoWiltonGlisson_ElectromagneticScatteringBySurfacesOfArbitraryShape} are applied to the \ac{MoM} formulation of the \ac{EFIE}~\cite{Gibson_MoMinElectromagnetics}, see red lines in Fig.~\ref{fig1}c. The algebraic form of the system equation reads
\begin{equation}
\Zm \Iv = \Vv,
\label{eq:MoM1}
\end{equation}
where $\Zm \in \mathbb{C}^{\Ndof \times \Ndof}$ is the system (impedance) matrix, $\Iv\in\mathbb{C}^{\Ndof\times 1}$ is a column vector of expansion coefficients for the current density (currents), and~$\Vv \in \mathbb{C}^{\Ndof\times 1}$ is a vector of excitation coefficients~\cite{Harrington_FieldComputationByMoM}. The inversion of \eqref{eq:MoM1}, \ie{},
\begin{equation}
\Iv = \Zm^{-1} \Vv = \Ym \Vv
\label{eq:MoM2}
\end{equation}
with $\Ym$ being the admittance matrix, gives the solution for the prescribed excitation. Formula~\eqref{eq:MoM2} binds the shape and the excitation with an induced current and its physical effects and, as such, is considered the key optimization constraint in this work.

Consider further that the excitation~$\Vv$ is given prior to the optimization and is fixed. Then, the only unknown in~\eqref{eq:MoM2} is the shape of the structure encoded in the impedance matrix~$\M{Z}$. Its parameterization is proposed in the next section.

For a particular state of the parametrization, there is a unique current vector~$\M{I}$ for a given excitation~$\Vv$. The optimized metric~$f \left( \Iv \right)$ or any $k$-th optimization constraint~$h_k \left( \Iv \right) = 0$ can then be defined based on this current vector. In the cases treated in this paper, functions~$f$ and~$h_k$, representing quadratic or linear forms in the current vector~$\M{I}$, are used to express quality factor, antenna gain, input impedance, or cycle mean powers.

\subsection{Generic Optimization Problem}
\label{sec:ShapeOpt:Opt}
\begin{figure}
\centering
\includegraphics[]{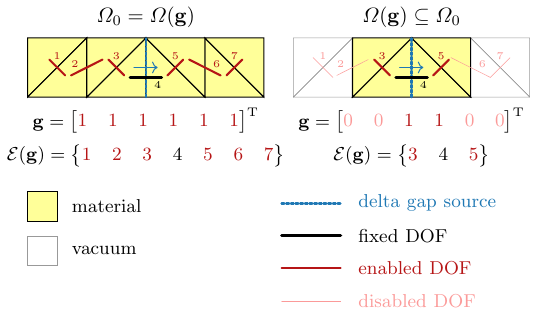}
\caption{Word~$\gene$ of an arbitrary shape within the design region~$\bbox$ expressed in a fixed discretization. (left) The entire structure with basis function number~4 as a fixed DOF. (right) Substructure formed by removing \ac{DOF} $\{1,2,6,7\}$ from the full strip.}
\label{fig3}
\end{figure}

Impedance matrix~$\Zm$ is of size $\Ndof \times \Ndof$, and the natural choice for a fixed discretization is to choose
\begin{equation}
\Nopt = \Ndof - P,
\label{eq:NA}
\end{equation}
optimization variables, where $P$ represents the number of \ac{DOF} which are considered to be fixed and whose state will not be optimized. The highest resolution of the optimization procedure is achieved when the same \ac{DOF} as for the impedance matrix representation, \ie{}, the basis functions~\cite{Capeketal_ShapeSynthesisBasedOnTopologySensitivity}, are chosen see Fig.~\ref{fig1}c. Consequently, any shape is represented by a finite-length \Quot{word} of $\Nopt$ symbols as
\begin{equation}
\gene = \begin{bmatrix}
g_1 & \cdots & g_n & \cdots & g_\Nopt
\end{bmatrix}^\trans,
\label{eq:param}
\end{equation}
where each \Quot{symbol}~$g_n$ represents material states of corresponding \ac{DOF},~$\gene$ fully describes a shape being optimized, see Fig.~\ref{fig3}, and the superscript ${}^\trans$ denotes matrix transpose. Notice that each state $g_n$ can take one of many values, for example, values indicating the presence of gold, silver, silica, vacuum, etc. For the sake of simplicity, and without loss of generality, each symbol in this paper is binary, $\gene \in \left\{0, 1\right\}^{\Nopt\times 1}$, indicating the presence of either a vacuum or a conductor. The word~$\gene$ fully populated by ones delimits the entire region used for optimization, called here design region and denoted~$\bbox$.

Having the optimized shapes properly represented, we can set up a general optimization problem~\cite{NocedalWright_NumericalOptimization} as
\begin{equation}
\begin{split}
\underset{\gene}{\T{minimize}} \quad & f \left( \Iv \left(\gene\right) \right) \\
\T{subject\, to} \quad & h_k \left( \Iv \left(\gene\right) \right) = 0,\quad k=1,...,K \\
& \Iv \left(\gene\right) = \Zm^{-1} \left(\gene\right) \Vv \\
& \gene \in \left\{0,1\right\}^{\Nopt \times 1},
\end{split}
\label{eq:optim1}
\end{equation}
where the only unknown is the binary word~$\gene$ indicating which \ac{DOF} are enabled (present) and which are disabled (removed) from the evaluation of~\eqref{eq:MoM2}, see Fig.~\ref{fig3}, hence describing the shape. Functions~$f$, $h_k$ depend on current~$\Iv\left(\gene\right)$ which is a subject of constraint~\eqref{eq:MoM2}. Consequently, current~$\Iv(\gene)$ flowing through the disabled edges ($g_b = 0$) is zeroed. The discretization is fixed for the entire optimization to lock the parameterization and to save computational time required for the recalculation of all the matrices whenever the meshing is changed.

\subsection{Combinatorial Topology Optimization}
\label{sec:ShapeOpt:OptB}
A notable property stemming from the formulation~\eqref{eq:optim1} is that the equality constraints cannot be fulfilled. There are $2^\Nopt$ possible values of $h_k \left(\Iv\left(\gene\right)\right)$ and, mathematically, there is no reason to expect that one of them is (exactly) equal to zero. However, we dispose of an extremely large variable space, so we can typically approach the required value closely in a finite-precision arithmetic~\cite{Sauer_NumericalAnalysis}.

In order to effectively overcome problems with multiple constraints and their fulfilment, the problem~\eqref{eq:optim1} is rewritten into a form compatible with discrete optimization by aggregating all metrics into a composite function~\cite{2008EichfelderAdaptiveScalarizationMethodsInMultiobjectiveOptimization}
\begin{equation}
\begin{split}
\underset{\gene}{\T{minimize}} \quad & f \left( \Iv \left(\gene\right) \right) + \sum_{k=1}^K w_k \left| h_k \left( \Iv \left(\gene\right) \right) \right| \\
\T{subject\, to} \quad & \Iv \left(\gene\right) = \Zm^{-1} \left(\gene\right) \Vv \\
& \gene \in \left\{0,1\right\}^{\Nopt \times 1},
\end{split}
\label{eq:optim2}
\end{equation}
where the weights~$w_k$ give additional freedom for the optimization, typically being swept to form a Pareto frontier~\cite{1978CohonMultiobjectiveProgrammingAndPlanning,Deb_MultiOOusingEA}. This trade-off also agrees well with what engineers do and expect -- the target criterion is counter-weighted by the fulfilment of constraints. Formulation~\eqref{eq:optim2} is used throughout this paper whenever topology optimization results are reported.

The combinatorial optimization problem~\eqref{eq:optim2} cannot be solved in polynomial time~\cite{Nemhauser_etal_IntegerAndCombinatorialOptimization} and only approximate solutions are available. Common strategies are heuristic binary searches using genetic algorithm (GA)~\cite{RahmatMichielssen_ElectromagneticOptimizationByGenetirAlgorithms} or density-based topology optimization~\cite{BendsoeSigmund_TopologyOptimization} with filtering process and subsequent thresholding. The strategy adopted in this paper uses a combination of a local algorithm based on discrete topology sensitivity and a genetic algorithm which preserves the diversity of the search space. This solution process is described in Sections~\ref{sec:Reanalysis},~\ref{sec:GA}~and~\ref{sec:memetics} and a major advantage of this approach is its discrete nature with no need for thresholding.

\subsection{Current Optimization and Fundamental Bounds}
\label{sec:ShapeOpt:FB}

An alternative approach to analyzing optimization problems~\eqref{eq:optim1} or~\eqref{eq:optim2} is to estimate their lower bounds which can also help to terminate the memetic routine discussed above. This leads to the concept of fundamental bounds, which already found its application in many antenna and scattering problems~\cite{GustafssonTayliEhrenborgEtAl_AntennaCurrentOptimizationUsingMatlabAndCVX, JelinekCapek_OptimalCurrentsOnArbitrarilyShapedSurfaces, GustafssonCapek_MaximumGainEffAreaAndDirectivity, 2020_Gustafsson_NJP}. 

In the framework of fundamental bounds, the optimization problem~\eqref{eq:optim1} is relaxed to determine its fundamental lower bound. The first step is a change of optimization variable from a word~$\gene$, which directly modifies the operators describing the physical problem, to current expansion coefficients~$\Iv$. The second step is a relaxation of the constraint~\eqref{eq:MoM1}, which only allows one solution for given material distribution and feeding. This constraint is relaxed to a complex power balance~$\Iv^\herm \Zm \Iv = \Iv^\herm \Vv$, which offers more freedom~\cite{2020_Gustafsson_NJP}. This finally leads to an optimization problem
\begin{equation}
\begin{split}
\underset{\Iv}{\T{minimize}} \quad & f \left( \Iv \right) \\
\T{subject\, to} \quad & h_k \left( \Iv \right) = 0, \quad k=1,...,K \\
& \Iv^\herm \Zm \Iv = \Iv^\herm \Vv,
\end{split}
\label{eq:optim1Bound}
\end{equation}
the solution of which can in many cases be found by methods of convex optimization~\cite{BoydVandenberghe_ConvexOptimization}, typically utilizing Lagrange duality~\cite{BoydVandenberghe_ConvexOptimization}.

It is important to stress here that any solution, $\Iv(\gene)$, to problem~\eqref{eq:optim1} is a solution to problem~\eqref{eq:optim1Bound} whenever problem~\eqref{eq:optim1} is solved on a material sub-region of that considered in~\eqref{eq:optim1Bound}. Due to its larger freedom, the problem~\eqref{eq:optim1Bound} thus always leads to a lower value of the objective function~$f$ than is allowed in~\eqref{eq:optim1}.

Unlike problems~\eqref{eq:optim1} or~\eqref{eq:optim2}, the solution to optimization problem~\eqref{eq:optim1Bound} is, when formulated as a dual convex problem, computationally inexpensive. The price to pay is a gap between the fundamental bound and solution to~\eqref{eq:optim2} which, nevertheless, is in many cases insignificant by means of practical design considerations, see Part~2 of this paper~\cite{2021_capeketal_TSGAmemetics_Part2}.

\subsection{An Illustrating Example}
\label{sec:ShapeOpt:Fitness}

To show the essential difference between the combinatorial problem~\eqref{eq:optim1} and its fundamental bound~\eqref{eq:optim1Bound}, a well-known example of minimizing the radiation Q-factor is shown in this section. The optimization problems for other metrics relevant to antenna design can be formulated analogously. 

The problem describing the minimal radiation Q-factor for an antenna in self-resonance reads~\cite{CapekGustafssonSchab_MinimizationOfAntennaQualityFactor}
\begin{equation}
\begin{split}
\underset{\gene}{\T{minimize}} \quad & Q_\T{U} (\gene) \\
\T{subject\, to} \quad & Q_\T{E} (\gene) = 0 \\
& \Iv \left(\gene\right) = \Zm^{-1} \left(\gene\right) \Vv \\
& \gene \in \left\{0,1\right\}^{\Nopt \times 1},
\end{split}
\label{eq:optim3}
\end{equation}
where
\begin{equation}
Q_\T{U} (\gene) = \dfrac{1}{2} \dfrac{\Iv^\herm \left(\gene\right) \Wm \Iv \left(\gene\right)}{\Iv^\herm \left(\gene\right) \Rmvac \Iv \left(\gene\right)}
\label{eq:QU}
\end{equation}
is an untuned Q-factor,
\begin{equation}
Q_\T{E} (\gene) = \dfrac{\left|\Iv^\herm \left(\gene\right) \Xm_0 \Iv \left(\gene\right)\right|}{\Iv^\herm \left(\gene\right) \Rmvac \Iv \left(\gene\right)}
\label{eq:QE}
\end{equation}
is the tuning part of the Q-factor~\cite{CapekJelinek_OptimalCompositionOfModalCurrentsQ}, and where $\Wm$ is the stored energy matrix~\cite{CapekGustafssonSchab_MinimizationOfAntennaQualityFactor}, and $\Rmvac$ and $\Xmvac$ are real and imaginary parts of the vacuum impedance matrix $\Zmvac = \Rmvac + \J \Xmvac$, respectively. The excitation vector~$\Vv$ typically corresponds to a localized source which, in this paper, is modeled by an excitation via a delta gap with voltage~$V_\T{in}$.

Following Section~\ref{sec:ShapeOpt:OptB}, the optimization problem~\eqref{eq:optim3} is rewritten as
\begin{equation}
\begin{split}
\underset{\gene}{\T{minimize}} \quad & Q (\gene) = w_1 Q_\T{U} (\gene) + w_2 Q_\T{E} (\gene) \\
\T{subject\, to} \quad & \Iv \left(\gene\right) = \Zm^{-1} \left(\gene\right) \Vv \\
& \gene \in \left\{0,1\right\}^{\Nopt \times 1},
\end{split}
\label{eq:optim4}
\end{equation}
where the choice of~$\left\{ w_1 = 1, w_2 = 0 \right\}$ yields an untuned Q-factor, while $\left\{w_1 = x_1, w_2 = 1-x_1\right\},\,0 < x_1 < 1$ yields a convex combination of a minimal untuned Q-factor and self-resonant behavior.

On the other hand, following Section~\ref{sec:ShapeOpt:FB}, and rewriting the power constraint in~\eqref{eq:optim1Bound} for a self-resonant ($Q_\T{E} = 0$ in~\eqref{eq:optim3}) localized source as
\begin{equation}
    \Iv^\herm \Zm \Iv = Y_\T{in}^* \left| V_\T{in} \right|^2 = 2P_{\T{r}},
\end{equation}
with~$Y_\T{in}$ denoting the input admittance seen by the source and~$P_{\T{r}}$ the radiated power, the fundamental bound on the lowest achievable Q-factor~\cite{CapekGustafssonSchab_MinimizationOfAntennaQualityFactor}, $Q_\T{lb}$, is found by
\begin{equation}
\begin{split}
\underset{\Iv}{\T{minimize}} \quad & \Iv^\herm \Wm \Iv \\
\T{subject\, to} \quad & \Iv^\herm \Rmvac \Iv = 2P_{\T{r}} \\
& \Iv^\herm \Xm \Iv = 0.
\end{split}
\label{eq:bnd1}
\end{equation}
This is a relaxation of~\eqref{eq:optim1Bound} valid for arbitrary feed positions and admittance values. Corresponding bounds, including constraints on the feed location and input admittance, are analogous by incorporating the constraints in $\M{V}$, \cf{}~\eqref{eq:optim1Bound}.

\begin{figure}
\centering
\includegraphics[]{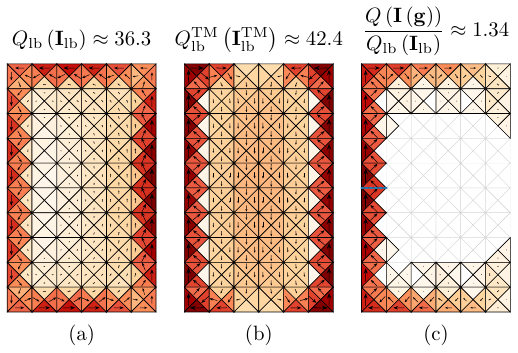}
\caption{Comparison of currents (colormap depicts the absolute value of surface current density, while arrows depict a snapshot) generated by the optimization problem for the lowest Q-factor found via the procedure from~\cite{CapekGustafssonSchab_MinimizationOfAntennaQualityFactor} (left), from~\cite{GustafssonTayliEhrenborgEtAl_AntennaCurrentOptimizationUsingMatlabAndCVX} (middle), and current obtained by an approximate solution to~\eqref{eq:optim4} with a delta gap feed highlighted by blue color (right). The electrical size of the rectangular region is~$ka = 0.5$ ($k$~is the wavenumber), the number of \ac{DOF} is $\Ndof=345$ (with one \ac{DOF} fixed as a feeding delta gap, \ie{}, $\Nopt = 344$). Q-factor for the lower bound (left) reaches~$Q_\T{lb} \approx 36.3$, for the TM-modes-constrained Q-factor~$Q_\T{lb}^\T{TM} \approx 42.4$, and $Q \approx 48.6$ for a shape found with the presented procedure of topology optimization, \ie{}, the performance of the topology optimized design is $1.34$ times above than the lower bound.}
\label{fig4}
\end{figure}

The optimal current density resulting from~\eqref{eq:bnd1} is shown in Fig.~\ref{fig4}a together with the corresponding Q-factor (lower bound on the Q-factor for the considered region). The middle panel of the same figure shows a result of a slightly altered problem~\cite{GustafssonTayliEhrenborgEtAl_AntennaCurrentOptimizationUsingMatlabAndCVX}, generating bound~$Q_\T{lb}^\T{TM}$, which is more representative for an electrically small single-port antenna since the radiation is restricted to TM spherical waves (electric dipoles). Finally, the current obtained by an approximate solution to~\eqref{eq:optim4} using the method described in this paper is shown in Fig.~\ref{fig4}c. 

Knowing the fundamental bound, the performance of a solution returned by an algorithm trying to solve problem~\eqref{eq:optim4} can be judged in an absolute sense by normalizing it to the fundamental bound
\begin{equation}
q (\gene) = \dfrac{Q(\gene)}{Q_\T{lb}}.
\label{eq:Qnorm}
\end{equation}
The termination criterion of the optimization algorithm is then defined as a threshold, $q(\gene) \leq t$, $t  \in \left[1, \infty\right)$, below which the performance is acceptable. Setting $t \rightarrow 1$ implies the actual shape~$\gene$ has to perform as well as the best hypothetical device. Such a requirement is typically not possible to achieve~\cite{Capek_etal_2019_OptimalPlanarElectricDipoleAntennas}.

\section{Gradient-based Local Step}
\label{sec:Reanalysis}

This section explains in detail the local step of a topology optimization based on a \ac{MoM} description and exact reanalysis. It is assumed that all \ac{MoM} matrices required for the evaluation of~\eqref{eq:optim1} are precalculated, the optimization problem is set up, and its lower bound is known. 

The technique introduced in this section is inspired by an exact reanalysis~\cite{Ohsaki_OptimizationOfFiniteDimensionalStructures} known in mechanics and used for fast evaluation of the smallest perturbations of the stiffness matrix~\cite{Jin_FEM2014}. It is an inversion-free technique, making it possible to investigate all the smallest topology changes in a computationally inexpensive way.

\begin{figure}
\centering
\includegraphics[]{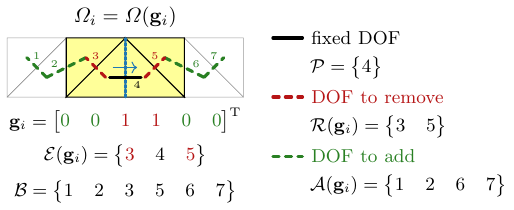}
\caption{Sketch of the sets enlisted in Table~\ref{Tab:TSnomenclature} and used for rank-1 inversion-free updates via exact reanalysis procedure exemplified for the shape in Fig.~\ref{fig3}. DOFs $\setRem(\gene_i)=\{3\ 5\}$ can be removed, $\setAdd(\gene_i)=\{1\ 2\ 6\ 7\}$ can be added, and $\setFixed=\{4\}$ is fixed.}
\label{fig7}
\end{figure}

The nomenclature is introduced first, building on Section~\ref{sec:ShapeOpt:Opt}. All \ac{DOF}, uniquely numbered from $n=1$ to $n=N$, see Fig.~\ref{fig3}, are for each shape~$\gene_i$ grouped into several sets, see Fig.~\ref{fig7}. Set~$\setGo$ contains all \ac{DOF} defining the design region (they are subjects of the optimization). Set~$\setFixed$ contains \ac{DOF} which are excluded from the optimization (feeding point, user-defined fixed part of the structure, etc.). Apart from sets~$\setGo$ and $\setFixed$, which are fixed for the entire duration of the optimization, there are sets~$\setGi(\gene_i)$, $\setAdd(\gene_i)$, and $\setRem(\gene_i)$, which are functions of the actual shape~$\gene_i$, where index $i$ denotes a particular state of the design. Set~$\setGi(\gene_i)$ collects all \ac{DOF}, which are currently part of the structure\footnote{In the case of multi-material optimization, each \ac{DOF} can represent one of many materials, it can, however, only represent one material at a time.}. Set~$\setAdd(\gene_i)$ contains all \ac{DOF} which might be added to the structure (\eg{}, to be changed from vacuum to material) and, conversely, set~$\setRem(\gene_i)$ contains all \ac{DOF} to be removed from the structure (\eg{}, to be changed from material to vacuum). This means, that there is~$\vert \setRem (\gene_i) \vert + \vert \setAdd (\gene_i) \vert$ possible shapes~$\gene_{i+1}$, which result from an update of shape~$\gene_i$ by a single \ac{DOF}. The meaning of all symbols is recapitulated in Table~\ref{Tab:TSnomenclature}.

\begin{table}[t] 
\centering
\caption{Sets used in this paper. The first two sets only depend on the design region and user's preferences, and the last three sets depend on the actual shape encoded by word~$\gene_i$}
\begin{tabular}{cll} 
set & description & note \\[1pt] \toprule
$\setGo$ & optimized \ac{DOF} & $\vert \setGo \vert = B$ \\
$\setFixed$ & \Quot{fixed} \ac{DOF} &  $\vert \setFixed \vert = P$ \\ \midrule
$\setGi(\gene_i)$ & \ac{DOF} of actual shape &  \\
$\setRem(\gene_i)$ & \ac{DOF} \Quot{to remove} & $\setRem(\gene_i) \equiv \setGi (\gene_i) \setminus \setFixed$ \\
$\setAdd(\gene_i)$ & \ac{DOF} \Quot{to add} & $\setAdd(\gene_i) \equiv \setGo \setminus \setGi (\gene_i)$ \\ \bottomrule
\end{tabular} 
\label{Tab:TSnomenclature}
\end{table}

A naive approach, conventionally followed by the pixeling technique~\cite{RahmatSamii_Kovitz_Rajagopalan-NatureInspiredOptimizationTechniques} based on heuristic optimization~\cite{RahmatMichielssen_ElectromagneticOptimizationByGenetirAlgorithms}, is to repetitively truncate impedance matrix~$\Zm$ and excitation vector~$\Vv$, \ie{}, to  evaluate
\begin{equation}
\Iv_\setGi = \Ym_\setGi \Vv_\setGi = \left(\Cm_\setGi^\trans \Zm \Cm_\setGi \right)^{-1} \Cm_\setGi^\trans \Vv,
\label{eq:pixelMoMinv}
\end{equation}
with an indexing matrix~$\Cm_\setGi$ truncating unused \ac{DOF}, so that various shapes represented by word~$\gene_i$ are iteratively studied, see Fig.~\ref{fig3}. Each trial requires matrix inversion, which has algorithmic complexity~$\OP{O}(M^3)$ where~$M$ is the dimension of the impedance matrix~$\Zm_\setGi = \Ym_\setGi^{-1}$ (number of nonzero entries in word~$\gene_i$). The cubic complexity of just one trial renders the entire optimization process computationally demanding and, in addition to this, heuristic updates have no information regarding sensitivity to local topology perturbations, so the convergence may be slow or none~\cite{Sigmund_OnTheUselessOfNongradinetApproachesInTopoOptim}. Both these deficiencies are eliminated with the procedure described below.

\subsection{Topology Sensitivity}
\label{sec:TopoSens}

\begin{figure}
\centering
\includegraphics[]{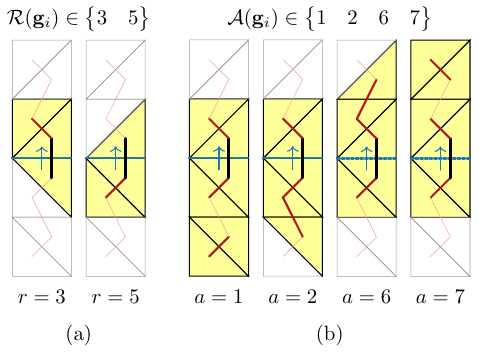}
\caption{All possible smallest perturbations~$\gene_{i+1}$ of shape~$\gene_i$ shown in Fig.~\ref{fig7} grouped into possible removals of the $r$-th \ac{DOF} from the set of all possible removals~$\setRem(\gene_i)$ (a), and possible additions of the $a$-th \ac{DOF} from the set of all possible additions~$\setAdd(\gene_i)$ (b). The $4$-th \ac{DOF} is considered to be fixed, $\setFixed \in \left\{4\right\}$, and its state cannot be changed.}
\label{fig5}
\end{figure}

\begin{figure}
\centering
\includegraphics[width=\columnwidth]{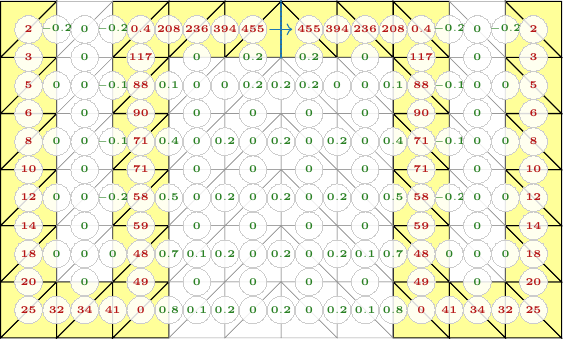}
\caption{Topology sensitivity map, $\V{\tau}(f,\gene_i)$, for an \textit{ad hoc} designed meandered dipole represented in discretization from Fig.~\ref{fig4}. Q-factor is~$Q(\Iv(\gene)) \approx 74.4$, the fundamental bound is~$Q_\T{lb} \approx 36.3$, see Fig.~\ref{fig4}. The potential addition of a \ac{DOF} is depicted in green, while the potential removal of a \ac{DOF} is depicted in red. Numerical values indicate an increase (positive sign) or a decrease (negative sign) in Q-factor associated with given action over \ac{DOF}. The depicted design is close to a local extremum allowing only a minuscule reduction of Q-factor by adding a \ac{DOF}.}
\label{fig6}
\end{figure}

For shape~$\gene_i$, as in Fig.~\ref{fig3}, we propose to iteratively evaluate all the smallest perturbations, \ie{}, to investigate the performance of all shapes~$\gene_{i+1}$ with Hamming distance~\cite{Robinson_AnIntroductionToAbstractAlgebra}~$d_\T{H}\left(\gene_i, \gene_{i+1}\right) = 1$, see Fig.~\ref{fig5}. The effects of the smallest topology perturbations over shape~$\gene_{i}$ are quantified as~\cite{Capeketal_ShapeSynthesisBasedOnTopologySensitivity}
\begin{equation}
\V{\tau} \left( f, \gene_i \right) = - \left[ f \left(\Iv\left(\gene_i\right), \gene_i \right) - f\left(\Iv \left(\gene_{i+1}\right), \gene_{i+1}\right) \right],
\label{eq:topoSens1}
\end{equation}
where the right-hand side is evaluated for all the smallest perturbations (see the following subsections), a process which is greatly accelerated by the employed low-rank updates. Vector~$\V{\tau}$ measures topology sensitivity of objective function~$f$ over shape~$\gene_i$, see Fig.~\ref{fig6}. It provides, on a particular case of Q-factor, information on what happens if a removed \ac{DOF} is put back (changed from vacuum to a material) or one active \ac{DOF} is removed (changed from material to vacuum). Objective function~$f$ can be a composite function and can consist of parts depending on current~$\Iv\left(\gene_i\right)$ or directly on word~$\gene_i$. The computational complexity of one local update (to evaluate one modified shape~$\gene_{i+1}$) is of the order~$\OP{O}(M)$ and computational complexity of all local updates (to investigate all possible modifications of shape~$\gene_i$) is of the order~$\OP{O}(M \Nopt)$, where $\Nopt = \Ndof - P$. In other words, gathering information about all the smallest perturbations~\eqref{eq:topoSens1} costs as much as a single solution to problem~\eqref{eq:MoM1}.

\subsection{Rank-1 Modifications and Updates}
\label{sec:TopoSens:rank1}

The smallest topology perturbations~\eqref{eq:topoSens1} are evaluated with block inversion~\cite{GolubVanLoan_MatrixComputations}
\begin{equation}
\begin{split}
\Ym & = \left[ \begin{array}{cc}
\Ym_{11} & \Ym_{12} \\
\Ym_{21} & \Ym_{22}
\end{array} \right]
=\left[ \begin{array}{cc}
\Zm_{11} & \Zm_{12} \\
\Zm_{21} & \Zm_{22}
\end{array} \right]^{-1} \\
& = \left[ \begin{array}{cc}
\Zm_{11}^{-1} + \Zm_{11}^{-1} \Zm_{12} S^{-1} \Zm_{21} \Zm_{11}^{-1} & - \Zm_{11}^{-1} \Zm_{12} S^{-1} \\
- S^{-1} \Zm_{21} \Zm_{11}^{-1} & S^{-1}
\end{array} \right]
\label{eq:SMW1}
\end{split}
\end{equation}
with the Schur complement~\cite{GolubVanLoan_MatrixComputations}
\begin{equation}
S = \Zm_{22} - \Zm_{21} \Zm_{11}^{-1} \Zm_{12},
\label{eq:SMW2}
\end{equation}
where~$\M{Z}_{mn}$ refers to a block of matrix~$\M{Z}$ and the upper-left block matrix of~\eqref{eq:SMW1} is called the Sherman-Morrison-Woodbury formula~\cite{1989_Hager_SIAM_R}.

For a \ac{DOF} addition, the identity is used to directly evaluate~$\Ym$ from~$\Zm_{11}^{-1}$ and matrix-vector products involving $\Zm_{11}^{-1}$. For a \ac{DOF} removal, the identity is used backwards to evaluate~$\Zm_{11}^{-1}$ from~$\Ym$ as the Schur complement $\Zm_{11}^{-1}=\Ym_{11}-\Ym_{12}\Ym_{22}^{-1}\Ym_{21}$. In both cases, only summations and matrix-vector products are needed.

The technical details covering all algebraic manipulations are summarized in Appendices~\ref{sec:TopoSens:RemPert}--\ref{sec:TopoSens:AddUpdate}. Appendices~\ref{sec:TopoSens:RemPert} and~\ref{sec:TopoSens:AddPert} describe removal of $r$-th and addition of $a$-th \ac{DOF}, respectively, \ie{}, how to get a currents~$\Iv_{\setGi\setminus r}$ and~$\Iv_{\setGi\cup a}$, representing perturbed shapes. Appendix~\ref{sec:TopoSens:Metrics} shows how to evaluate optimized metrics based on current matrices~$\left[ \Iv_\setRem (\gene_i) \right]$ and $\left[ \Iv_\setAdd (\gene_i) \right]$, where all possible perturbations of a word~$\gene_i$ were collected. Finally, Appendices~\ref{sec:TopoSens:RemUpdate} and~\ref{sec:TopoSens:AddUpdate} summarize how to update the shape once the locally-optimal perturbation is found. The local optimization algorithm is visually depicted in Fig.~\ref{fig8}. The code with all formulas implemented is available in supplementary material~\cite{TSGA22_SuplMat}.

\begin{figure}
\centering
\includegraphics[width=\columnwidth]{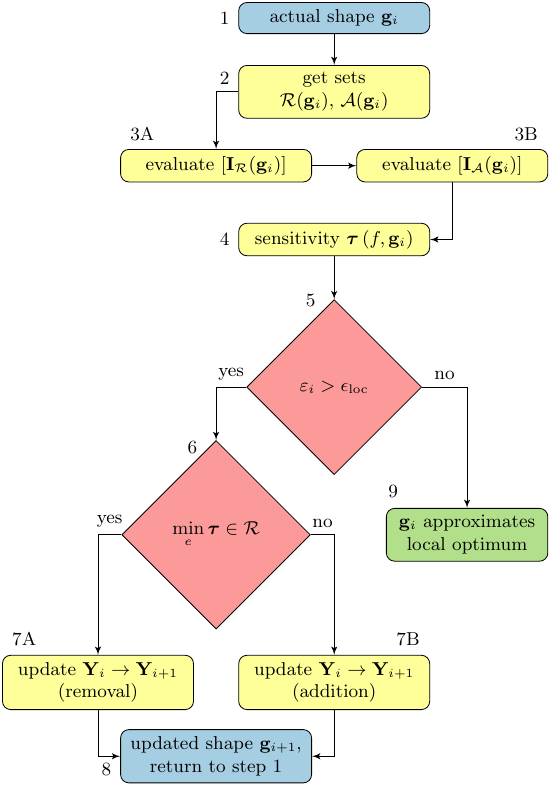}
\caption{Flowchart of topology sensitivity evaluation via local perturbations and a local update of a shape. An initial shape is represented by~$\gene_i$ (step~1). All possible removals and additions are collected in sets~$\setRem(\gene_i)$ and~$\setAdd(\gene_i)$ (step~2). Topology sensitivity~$\M{\tau}(f,\gene_i)$ of an optimized metric to geometry perturbations is evaluated (step~4) from currents~$\left[\Iv_\setRem(\gene_i)\right]$ and~$\left[\Iv_\setAdd(\gene_i)\right]$ (step~3). During the local optimization, the shape is iteratively modified by returning from step~8 to step~1 as long as there are negative sensitivities available (as long as step~9 is reached). The local step can, alternatively, be terminated sooner, based on additional terminal criteria.}
\label{fig8}
\end{figure}

\section{Memetic Algorithm For Topology Optimization}
\label{sec:memetics}

\begin{figure}
\centering
\includegraphics[]{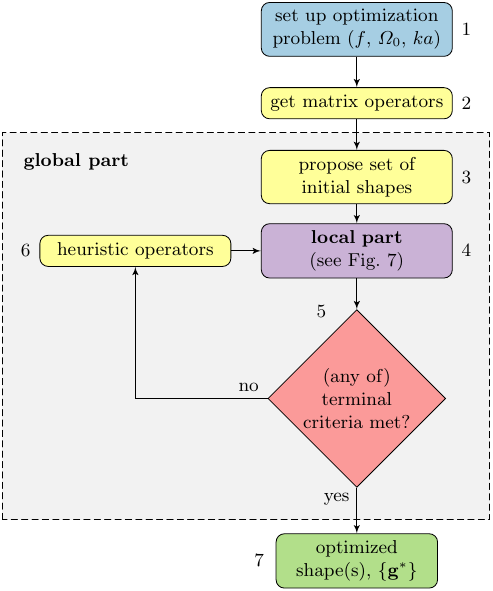}
\caption{Schematics of the memetic topology optimization algorithm based on exact reanalysis of \ac{MoM} models. After setting up the optimization problem (step~1), the matrix operators required to solve the problem are evaluated (step~2). The optimization starts with an initial set of seeds (step~3) which are locally optimized (step~4) and further improved with heuristic operators (step~6). The entire procedure is repeated until the terminal criteria are met (step~5).}
\label{fig9}
\end{figure}

The topology sensitivity introduced in Section~\ref{sec:Reanalysis} can be utilized to iteratively update the geometry represented by word~$\gene_i$ until the position of local minimum is found with sufficient precision. The local optimization procedure from Fig.~\ref{fig8} is repeated as long as the relative improvement of the objective function, defined as
\begin{equation}
    \varepsilon_i = -\dfrac{1}{\vert f_\T{min} \vert} \min \left\{ \V{\tau} \left( f, \gene_i \right) \right\},
\end{equation}
where $f_\T{min}$ is the minimum value of the objective function from the $i$-th iteration, is higher than a predefined limit~$\epsilon_\T{loc}$, \ie{}, as long as $\varepsilon_i > \epsilon_\T{loc}$. The problem of inverse design is, however, non-convex. To make the procedure robust and to prevent it from getting stuck in a local minimum, a global optimization step is considered as well, resulting in a memetic optimization scheme, see~Fig.~\ref{fig9}.

The memetic optimization procedure~\cite{Hart_etal_RecentAdvancesInMemeticAlgorithms} combines two approaches to increase the robustness and convergence of the optimization. Since the type of the problem to be solved is unknown, it is advantageous to combine as different algorithms as possible~\cite{WolpertMacready_NoFreeLunchTheoremsForOptimization}. In this work, it is the gradient-based topology sensitivity utilizing exact reanalysis of \ac{MoM} models~\cite{Ohsaki_OptimizationOfFiniteDimensionalStructures}, see Section~\ref{sec:Reanalysis}, and the genetic algorithm~\cite{RahmatMichielssen_ElectromagneticOptimizationByGenetirAlgorithms, Haupt_Werner_GeneticAlgorithmsInEM}.

\subsection{Heuristic Global Step}
\label{sec:GA}

Since word~$\gene_i$ is discrete and, in the case of this work, binary, a genetic algorithm (GA)~\cite{RahmatMichielssen_ElectromagneticOptimizationByGenetirAlgorithms} is selected as a global scheme. The algorithm implemented in a MATLAB package FOPS~\cite{fops} is utilized, enabling direct communication with the local optimization routines implemented in the same environment. 

Initial generation~$\mathcal{G}_1$ (a set of words~$\gene$ for the first global iteration, $j = 1$) with size~$\Nags$ is produced randomly within the decision space, see step 3 in Fig.~\ref{fig9}. Two extreme cases having all letters in word~$\gene$ being equal to vacuum and metal are always included in this initial set as they represent physically meaningful states (there is either no device or a device spanning the entire design region). The global iterative process consists of steps~4-6 in Fig.~\ref{fig9}. In every iteration~$j$, a new generation~$\mathcal{G}_j$ is created using three fundamental GA-based heuristic operations~\cite{1997_Johnson_Genetic}: mating pool selection, crossover, and mutation, see step~6 in Fig.~\ref{fig9}. 

We prefer diversity rather than global convergence in GA because the local step in Section~\ref{sec:Reanalysis} provides much faster overall convergence. Therefore, unlike conventional approaches~\cite{JohnsonSamii1997_GAinEM}, all the heuristic operations are set to maximize the diversity of the new population.  The binary tournament selection~\cite{1997_Johnson_Genetic} was implemented to select the mating pool. 
The crossover operation is then applied to the words selected from the mating pool to provide a new set~$\mathcal{G}_{j+1}$. A probability of crossover is set in this study as~$p_{\mathrm{C}} = 0.9$ to enhance further the diversity among words in the new set $\mathcal{G}_{j+1}$. The operation mutation prevents the GA from getting stuck in a local minimum. The operation takes a word from new generation~$\mathcal{G}_{j+1}$, selects a randomly chosen letter, and flips its value. The mutation is applied with a probability~$p_{\mathrm{M}}$. The maximal value~$p_{\mathrm{M}} = 1$ was used in this study to enhance the diversity among present words. This is crucial, especially in later iterations of the algorithm, when most words become very similar or the same. This behavior results from the local method where a larger number of initial words are modified into the same word by the local algorithm, \ie{}, the words are attracted into the same local minimum.

\subsection{Initiation}
\label{sec:TopoSens:Start}

Before optimization begins, region~$\srcRegion$ is discretized with a required granularity, see Fig.~\ref{fig1}. The number of~\ac{DOF}~$\Ndof$ and the number of optimization unknowns~$\Nopt$ predetermine the resolution of the optimization. Depending on the implementation and type of the problem, the algorithmic complexity is between~$\OP{O}\left(\Ndof^3\right)$ and~$\OP{O}\left(\Ndof^5\right)$. Increasing~$\Ndof$, however, also enlarges the solution space and ensures a better performance from the optimized shape.

The next step is to evaluate all required \ac{MoM} operator matrices, and to define an excitation. Then, the objective function is formulated. When a composite function is chosen, the weighting coefficients are selected.

The optimization has many control parameters, see Table~\ref{tab:controlParam} for an overview. These have to be set at the beginning and strongly influence the course of the optimization and its results. For example, it is not necessary to run the local step (greedy search) until the local minimum is found. Instead, the maximum number of iteration~$I$ is specified, or, alternatively, the local updates are performed only until the relative difference between values of the objective function evaluated between two consecutive iterations~$i$ and $i+1$ is smaller than a predefined value~$\epsilon_\T{loc}$. Another possibility is to disable either removals or additions by imposing~$\setRem = \emptyset$ or $\setAdd = \emptyset$, respectively. The maximum number of iterations and the relative difference can be specified for the global step as well. In this case, relative difference~$\epsilon_\T{glob}$ is evaluated for the worst performing agent of two consecutive iterations~$j$ and~$j+1$. An important parameter is the number of agents~$\Nags$. The probabilities of crossover and mutation is fixed to~$1$ in this work, but can be set to another number, or can be changed during the optimization.

\begin{table}[]
\caption{Inputs and control parameters for memetic algorithm}
\centering
\begin{tabular}{ccc} 
parameter & description \\ [0.5ex] 
\toprule
\bf{local step} & (Greedy search based on rank-1 topology differences) \\ \midrule
$I$ & maximum number of local iterations \\ [0.5ex]
$\epsilon_\T{loc}$ & relative difference to terminate local step \\ [0.5ex]
$\setRem \equiv \emptyset$ & removals are disabled \\ [0.5ex]
$\setAdd \equiv \emptyset$ & additions are disabled \\ \toprule
\bf{global step} & (heuristic genetic algorithm) \\ \midrule
$J$ & maximum number of global iterations \\ [0.5ex]
$\epsilon_\T{glob}$ & relative difference to terminate global step \\ [0.5ex]
$\Nags$ & number of agents \\ [0.5ex]
$p_\T{C}$ & crossover probability ($p_\T{C} = 1$) \\ [0.5ex]
$p_\T{M}$ & mutation probability ($p_\T{M} = 1$) \\ \toprule
\bf{memetics} & (combination of local and global steps) \\ \midrule
$c_\T{bnd}$ & acceptable distance from fundamental bound \\ [0.5ex]
$f$ & (composite) objective function \\ [0.5ex]
$\left\{ w_k, w_l \right\}$ & weights (if objective function if composite) \\ \bottomrule
\end{tabular}
\label{tab:controlParam}
\end{table}

\subsection{Termination}
\label{sec:TopoSens:TermCrit}

There are two sets of possible terminal criteria, one for the local step and one for the global step, see Table~\ref{tab:termParam}. The ultimate criterion for the local step is when the local minimum is found, \ie{}, all topology sensitivities are non-negative, $\tau_e < 0$ for any~$e$-th tested \ac{DOF}. However, it often happens that a large number of local updates improves the value of the objective function negligibly, just by tuning the fine details of the structure. This behavior is not required at the beginning of the optimization when the global step performs major changes of the structure on the global level. Therefore, to reduce computational burden, relative difference~$\epsilon_\T{loc}$ and the maximum number of iterations~$I$ should routinely be used. The same applies to the global step with~$\epsilon_\T{glob}$ and~$J$ parameters. As the final terminal criterion, the distance to the fundamental bound~$c_\T{bnd}$ is specified at the beginning. It has to be noted, that such performance is not achievable, and the optimization stops because of reaching either iteration~$J$ or relative difference~$\epsilon_\T{glob}$.

\begin{table}[]
\label{tab:termParam}
\caption{Possible terminal criteria for local and global steps}
\centering
\begin{tabular}{ccc} 
parameter & condition \\ [0.5ex] 
\toprule
\bf{local step} & \\ \midrule
local minimum reached & $\tau_e \geq 0 \quad \forall e$ \\[2ex]
relative difference & $\dfrac{f_i - f_{i+1}}{f_i} < \epsilon_\T{loc}$ \\[2ex]
number of iterations & $i > I$ & \\ \toprule
\bf{global step} & \\ \midrule
distance from fund. bound & $\dfrac{f_{i+1}}{f_\T{bound}} < c_\T{bnd}$ \\[3ex]
relative difference & $\dfrac{\max\{\V{f}_j\} - \max\{\V{f}_{j+1}\}}{\max\{\V{f}_i\}} < \epsilon_\T{glob}$ \\[3ex]
number of iterations & $j > J$ \\ \bottomrule
\end{tabular}
\end{table}

\subsection{Effective Implementation}
\label{sec:TopoSens:Implement}

The entire local step of one topology perturbation is depicted in Fig.~\ref{fig8} and is based on the evaluation of~\eqref{eq:rem-I}--\eqref{eq:add-V}. For the cases treated in this paper these formulas are typically evaluated billions of times per optimization run, so it is critical to implement them as effectively as possible:
\begin{itemize}
    \item All operations in \eqref{eq:rem-I}--\eqref{eq:add-V} are vectorized, \ie{}, all formulas are evaluated in a cycle-free procedure for all possible perturbations (removals and additions). Notice that the removals are significantly faster than additions, since each addition contains one matrix-vector product. When properly implemented, one entire evaluation of topology sensitivity should have a computational cost similar to a solution to~\eqref{eq:MoM2}. In other words, the evaluation of one \ac{MoM} problem is equivalent to checking all the smallest perturbations.
    \item All multiplications involving indexing matrices can be transformed into direct indexing, \ie{}, explicit matrix multiplication is not needed.
    \item The evaluation of the fitness function, which creates the biggest computational burden (contains matrix-matrix multiplication), might be evaluated on a GPU card.
\end{itemize}

The resulting, locally optimal, word is the input to the global step evaluation, see step~4 in Fig.~\ref{fig9}. The global algorithm, described in Section~\ref{sec:GA}, is performed multiple times with various starting words. The underlying local step has high granularity which favors the utilization of parallel computing.

\section{Conclusion}
\label{sec:concl}

A memetic optimization combining local and global steps was introduced. It solves the original combinatorial problem of the optimal material distribution in a prescribed region. In this paper, the smallest topological updates are inversion-free and derived for a method-of-moments system (impedance) matrix and related matrices. Thanks to the investigation of all the smallest perturbations in each step, the method is suitable for small and medium-sized mesh grids, \ie{}, for problems typically solved with a direct solver. This involves the majority of studies where small changes in geometry cause large changes in performance. There is no restriction on a particular form of an objective function. The computational cost of one iteration of the local step gathering all topology sensitivities is equivalent to one solution to the problem. The method converges fast due to the use of the local step. The global step is used to maintain diversity during the optimization. It operates in a greatly reduced solution space containing only locally optimal solutions effectively found by the local step. A unique advantage of this approach is the direct comparison of the achieved result with fundamental bounds.

\appendices
\section{Rank-1 Perturbation -- \ac{DOF} Removal}
\label{sec:TopoSens:RemPert}

The effect of a rank-1 \ac{DOF} removal is considered in this Appendix. Applying~\eqref{eq:SMW1} as shown in~\cite{Capeketal_ShapeSynthesisBasedOnTopologySensitivity}, and taking into account arbitrary excitation represented by excitation vector~$\Vv$, we get
\begin{equation}
\Iv_{\setGi\setminus r} = \Cm^\trans_{\setGi\setminus r} \left( \Iv_\setGi - \dfrac{I_r}{Y_{rr}} \M{y}_{\setGi\setminus r} \right) \Cm_{\setGi\setminus r},
\label{eq:rem-I}
\end{equation}
where subindex~$\setGi\setminus r$ denotes the removal of the $r$-th edge from the actual structure represented by a set of edges~$\setGi(\gene_i)$, see Table~\ref{Tab:TSnomenclature}, $I_r$ is a current flowing through the $r$-th basis function prior to the modification, $Y_{rr}$ is the $r$-th diagonal term of matrix~$\Ym_\setGi$, $\M{y}_{\setGi\setminus r}$ is the $r$-th column of matrix~$\Ym_\setGi$ and
$\Cm_{\OP{E} \setminus r} = [C_{\OP{E} \setminus r, nn}]$ is an indexing matrix,
\begin{equation}
\label{eq:topoRemCmat}
C_{\setGi\setminus r, nn} = \left\{
\begin{array}{lll}
0 & \Leftrightarrow & n = r, \\
1 & \Leftrightarrow & \mathrm{otherwise}, \\
\end{array}
\right.
\end{equation}
in which all columns full of zeros are removed. Current~$\Iv_\setGi$ corresponding to word~$\gene_i$ is supposed to be known from the previous iteration. When all possible removals~\eqref{eq:rem-I} are evaluated, the perturbed currents are collected in a matrix as
\begin{equation}
\left[ \Iv_\setRem (\gene_i) \right] = \left[ \Iv_{\setGi\setminus r} \right], \quad r \in \setRem(\gene_i).
\label{eq:rem-IR}
\end{equation}

A physical interpretation of~\eqref{eq:rem-I} is such that the original current~$\Iv_\setGi$ is zeroed at the $r$-th position by a placement of a particular voltage source. This might be understood as an application of Norton's equivalence theorem~\cite{Wing_ClassicalCircuitTheory}, zeroing the $r$-th column and row in admittance matrix~$\Ym_\setGi$, therefore, as a removal of the $r$-th column\footnote{Index~$r$, and similarly for index~$a$ introduced in Appendix~\ref{sec:TopoSens:AddPert}, has to be understood as relative indices. The matrices change their size during removals and additions and a current flowing through one fixed edge might change its position in vectors and matrices denoted by subindex~$\setGi$.} and $r$-th row of impedance matrix~$\Zm_\setGi$.

The implementation of \ac{DOF} removal is shown in function \textsc{removeEdges.m}, see Supplementary material~\cite{TSGA22_SuplMat}.

\section{Rank-1 Perturbation -- \ac{DOF} Addition}
\label{sec:TopoSens:AddPert}

Similar to removals introduced in Appendix~\ref{sec:TopoSens:RemPert}, the previously removed \ac{DOF} can be added back to a structure. For this purpose, the original impedance matrix~$\Zm$ and excitation vector~$\Vv$ have to be stored in the memory.

The current with the $a$-th \ac{DOF} added is evaluated as~\cite{Capeketal_InversionFreeEvaluationOfNearestNeighborsInMoM}
\begin{equation}
\Iv_{\setGi\cup a} = \Cm^\trans_{\setGi\cup a} \left( 
\left[ \begin{array}{*{10}{c}}
\Iv_\setGi \\
0
\end{array} \right] - 
\dfrac{v_a - \left(\Cm_\setGi^\trans \M{z}_a \right)^\trans \Iv_\setGi}{z_a}
\left[ \begin{array}{*{10}{c}}
\M{x}_a \\
0
\end{array} \right]
\right) \Cm_{\setGi\cup a}
\label{eq:add-I}
\end{equation}
where~$v_a$ is the voltage at the $a$-th position of the original vector~$\Vv$, the auxiliary variables~$\M{x}_a$ and~$z_a$ read
\begin{equation}
\M{x}_a = \Ym_\setGi \Cm_\setGi^\trans \M{z}_{a}
\quad\text{and}\quad
z_a = Z_{aa} - \left(\Cm_\setGi^\trans \M{z}_a \right)^\trans \M{x}_a,
\label{eq:add-Iaux1}
\end{equation}
respectively, with~$\M{z}_a$ being the $a$-th column of the original impedance matrix~$\Zm$, $Z_{aa}$ being the $a$-th diagonal term of the original matrix~$\Zm$, and the entries of indexing matrix~$\Cm_{\setGi \cup a} = [C_{\setGi \cup a, mn}]$ read
\begin{equation}
\label{eq:topoAddCmat}
C_{\setGi \cup a, mn} = \left\{
\begin{array}{lll}
1 & \Leftrightarrow & n = S\left(m\right), \\
0 & \Leftrightarrow & \mathrm{otherwise}, \\
\end{array}
\right.
\end{equation}
with $m \in \left\{ 1, 2, \dots, E+1\right\}$ and where $S$ is a set of target indices if a set~$\left\{ \OP{E}, a\right\}$ is sorted in ascending order. Finally, the perturbed currents for all possible additions are evaluated and collected as
$\left[ \Iv_\setAdd (\gene_i) \right] = \left[ \Iv_{\setGi\cup a} \right], \quad a \in \setAdd(\gene_i)$.

An implementation of the \ac{DOF} addition is shown in function \textsc{topoRemove.m}, see Supplementary material~\cite{TSGA22_SuplMat}.

\section{Evaluation of Optimization Metrics}
\label{sec:TopoSens:Metrics}

Once a structure is modified, the performance of a device is evaluated and the topology sensitivities~\eqref{eq:topoSens1} for all topology perturbations represented by currents~$\left[ \Iv_\setRem (\gene_i) \right]$ and~$\left[ \Iv_\setAdd (\gene_i) \right]$.

Since the dimension of matrices~$\left[ \Iv_\setRem (\gene_i) \right]$ and~$\left[ \Iv_\setAdd (\gene_i) \right]$ are different and the ordering of \ac{DOF} is different for each column, the matrix operators used for the evaluation of the objective function have to be modified according to
$\M{A}_e = \Cm_e^\trans \M{A} \Cm_e$,
where matrix~$\Cm_e$ coincides with~$\Cm_{\setGi\setminus r}$ for~$e = r$ (removal) and~$\Cm_{\setGi \cup a}$ for $e = a$ (addition) and cuts unused columns and rows from matrix~$\M{A}$, providing also correct ordering for the multiplication between the matrix operator and current vectors, \cf{}~\eqref{eq:optim1}. The algebraic operation $\Cm_e^\trans \M{A} \Cm_e$ is computationally inexpensive and is equivalent to direct indexing.

An implementation of the \ac{DOF} addition is shown in function \textsc{topoAdd.m}, see Supplementary material~\cite{TSGA22_SuplMat}.

\section{Rank-1 Update -- \ac{DOF} Removal}
\label{sec:TopoSens:RemUpdate}

The \ac{DOF} reaching the lowest (negative) value of sensitivity~$\V{\tau} \left( f, \gene_i \right)$, \ie{}, the highest decrease in objective function~$f$, \cf{}~\eqref{eq:topoSens1}, is to be updated at the end of the local sensitivity optimization step. The update is either a \ac{DOF} removal or a \ac{DOF} addition.

For the \ac{DOF} removal, admittance matrix~$\Ym_\setGi$ valid in $i$-th iteration is reduced as~\cite{Capeketal_ShapeSynthesisBasedOnTopologySensitivity}
\begin{equation}
\Ym_{\setGi\setminus r} = \Cm^\trans_{\setGi\setminus r} \left( \Ym_\setGi - \dfrac{1}{Y_{rr}} \M{y}_{\setGi\setminus r} \M{y}_{\setGi\setminus r}^\trans \right) \Cm_{\setGi\setminus r}
\label{eq:rem-Y}
\end{equation}
and the updated matrix~$\Ym_{\setGi\setminus r}$ is used for the next $i+1$ iteration. Consequently, the excitation vector is modified as well
$\Vv_{\setGi\setminus r} = \Cm^\trans_{\setGi\setminus r} \Vv_\setGi$.

An implementation of the system update by \ac{DOF} removal is shown in function \textsc{removeEdges.m}, see Supplementary material~\cite{TSGA22_SuplMat}.

\section{Rank-1 Update -- \ac{DOF} Addition}
\label{sec:TopoSens:AddUpdate}

For the \ac{DOF} addition, admittance matrix~$\Ym_\setGi$ valid in $i$-th iteration is expanded as~\cite{Capeketal_InversionFreeEvaluationOfNearestNeighborsInMoM}
\begin{equation}
\Ym_{\setGi\cup a} = \dfrac{1}{z_a} \Cm^\trans_{\setGi\cup a} 
\left[ \begin{array}{cc}
z_a \Ym_\setGi + \M{x}_a \M{x}_a^\trans & -\M{x}_a \\
-\M{x}_a^\trans & 1
\end{array} \right]
\Cm_{\setGi\cup a},
\label{eq:add-Y}
\end{equation}
and the excitation vector as
\begin{equation}
\Vv_{\setGi\cup a} = \Cm^\trans_{\setGi\cup a} \left[ \begin{array}{*{10}{c}}
\Vv_\setGi \\
V_a
\end{array} \right].
\label{eq:add-V}
\end{equation}

An implementation of the system update by \ac{DOF} addition is shown in function \textsc{addEdges.m}, see Supplementary material~\cite{TSGA22_SuplMat}.

\section*{Acknowledgement}
The access to the computational infrastructure of the OP VVV funded project CZ.02.1.01/0.0/0.0/16\_019/0000765 ``Research Center for Informatics'' is gratefully acknowledged.


\begin{IEEEbiography}[{\includegraphics[width=1in,height=1.25in,clip,keepaspectratio]{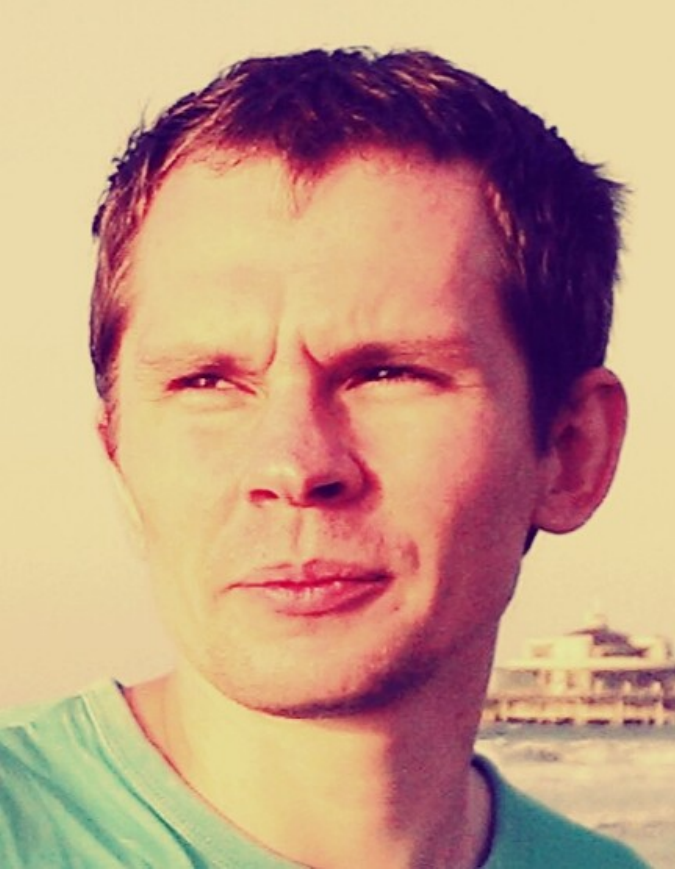}}]{Miloslav Capek}
(M'14, SM'17) received the M.Sc. degree in Electrical Engineering 2009, the Ph.D. degree in 2014, and was appointed Associate Professor in 2017, all from the Czech Technical University in Prague, Czech Republic.
	
He leads the development of the AToM (Antenna Toolbox for Matlab) package. His research interests are in the area of electromagnetic theory, electrically small antennas, numerical techniques, fractal geometry, and optimization. He authored or co-authored over 120~journal and conference papers.

Dr. Capek is Associate Editor of IET Microwaves, Antennas \& Propagation. He was a regional delegate of EurAAP between 2015 and 2020. He received the IEEE Antennas and Propagation Edward E. Altshuler Prize Paper Award 2023.
\end{IEEEbiography}

\begin{IEEEbiography}[{\includegraphics[width=1in,height=1.25in,clip,keepaspectratio]{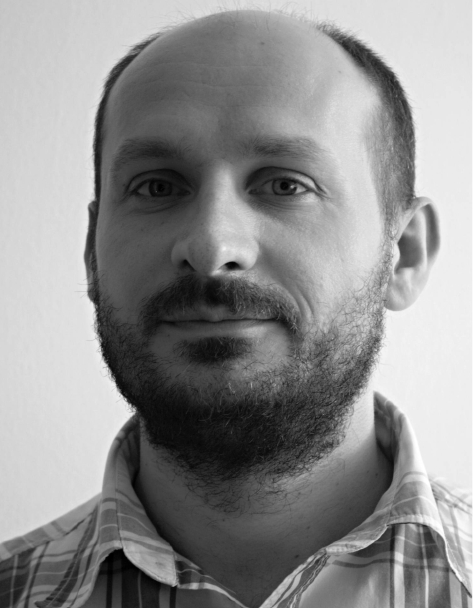}}]{Lukas Jelinek}
received his Ph.D. degree from the Czech Technical University in Prague, Czech Republic, in 2006. In 2015 he was appointed Associate Professor at the Department of Electromagnetic Field at the same university.

His research interests include wave propagation in complex media, electromagnetic field theory, metamaterials, numerical techniques, and optimization.
\end{IEEEbiography}

\begin{IEEEbiography}[{\includegraphics[width=1in,height=1.25in,clip,keepaspectratio]{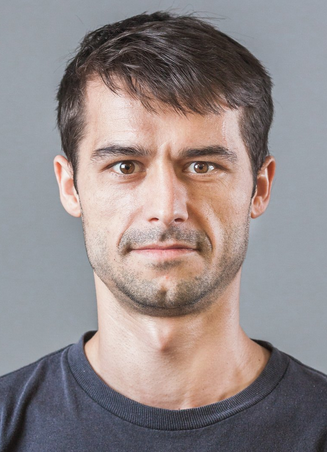}}]{Petr Kadlec}
(M'13) received the Ph.D. degree in electrical engineering from the Brno University of Technology (BUT), Brno, Czech Republic, in 2012. He is currently an Associate Professor with the Department of Radioelectronics, BUT. His research interests include global optimization methods and computational methods in electromagnetics. He is a leading developer of the FOPS (Fast Optimization ProcedureS) MATLAB software package.
\end{IEEEbiography}

\begin{IEEEbiography}[{\includegraphics[width=1in,height=1.25in,clip,keepaspectratio]{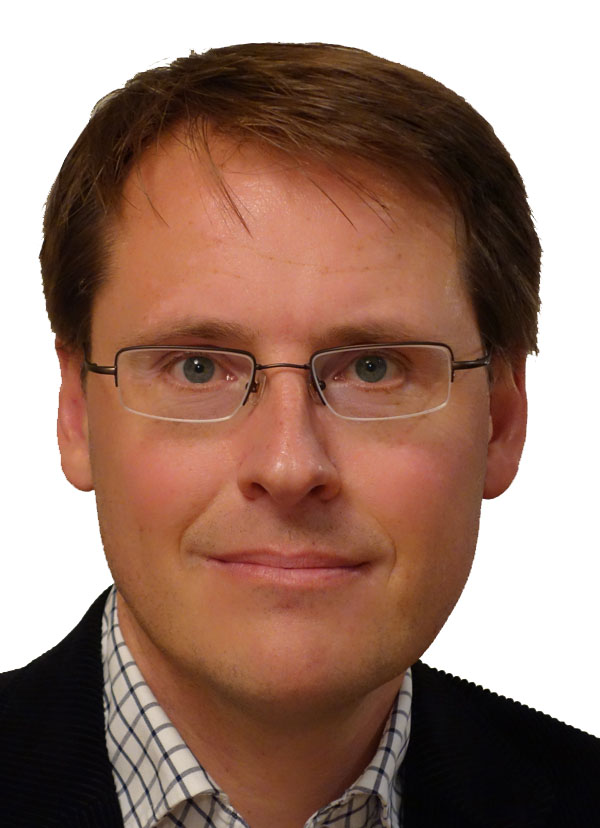}}]{Mats Gustafsson}
received the M.Sc. degree in Engineering Physics 1994, the Ph.D. degree in Electromagnetic Theory 2000, was appointed Docent 2005, and Professor of Electromagnetic Theory 2011, all from Lund University, Sweden.

He co-founded the company Phase holographic imaging AB in 2004. His research interests are in scattering and antenna theory and inverse scattering and imaging. He has written over 100 peer reviewed journal papers and over 100 conference papers. Prof. Gustafsson received the IEEE Schelkunoff Transactions Prize Paper Award 2010, the IEEE Uslenghi Letters Prize Paper Award 2019, and best paper awards at EuCAP 2007 and 2013. He served as an IEEE AP-S Distinguished Lecturer for 2013-15.
\end{IEEEbiography}

\end{document}